\documentclass{elsart}
\usepackage{amsmath}
\usepackage{amssymb}

\newcommand{\med}{{\rm median}}

\newcommand{\sug}{\mathcal{S}_{\mu}}

\begin{document}

\begin{frontmatter}

\title{Weighted lattice polynomials}

\author{Jean-Luc Marichal}
\ead{jean-luc.marichal[at]uni.lu}

\address{
Institute of Mathematics, University of Luxembourg\\
162A, avenue de la Fa\"{\i}encerie, L-1511 Luxembourg, Luxembourg }

\date{January 11, 2008}

\begin{abstract}
We define the concept of weighted lattice polynomial functions as lattice polynomial functions constructed from both variables and parameters.
We provide equivalent forms of these functions in an arbitrary bounded distributive lattice. We also show that these functions include the class
of discrete Sugeno integrals and that they are characterized by a median based decomposition formula.
\end{abstract}

\end{frontmatter}

\section{Introduction}

In lattice theory, \emph{lattice polynomials}\/ have been defined as well-formed expressions involving variables linked by the lattice
operations $\wedge$ and $\vee$ in an arbitrary combination of parentheses; see for instance Birkhoff \cite[\S II.5]{Bir67} and Gr\"atzer
\cite[\S I.4]{Grae03}. In turn, such expressions naturally define \emph{lattice polynomial functions}. For example,
$$
p(x_1,x_2,x_3)=(x_1\wedge x_2)\vee x_3
$$
is a 3-ary (ternary) lattice polynomial function.

The concept of lattice polynomial function can be straightforwardly generalized by fixing some variables as ``parameters'', like in the 2-ary
(binary) polynomial function
$$
p(x_1,x_2)=(c\vee x_1)\wedge x_2,
$$
where $c$ is a constant element of the underlying lattice.

In this paper we investigate those ``parameterized'' polynomial functions, which we shall call \emph{weighted lattice polynomial}\/ (w.l.p.)
functions. More precisely, we show that, in any bounded distributive lattice, w.l.p.\ functions can be expressed in disjunctive and conjunctive
normal forms and we further investigate these forms in the special case when the lattice is totally ordered. We also show that w.l.p.\ functions
include the discrete Sugeno integral \cite{Sug74}, which has been extensively studied and used in the setting of nonlinear aggregation and
integration. Finally, we prove that w.l.p.\ functions can be characterized by means of a median based system of functional equations.

Throughout, we let $L$ denote an arbitrary bounded distributive lattice with lattice operations $\wedge$ and $\vee$. We denote respectively by
$0$ and $1$ the bottom and top elements of $L$. For any integer $n\geqslant 1$, we set $[n]:=\{1,\ldots,n\}$ and, for any $S\subseteq [n]$, we
denote by $\mathbf{e}_S$ the characteristic vector of $S$ in $\{0,1\}^n$, that is, the $n$-dimensional vector whose $i$th component is $1$, if
$i\in S$, and $0$, otherwise. Finally, since $L$ is bounded,
$$\bigvee_{x\in\varnothing}x=0\quad\mbox{and}\quad\bigwedge_{x\in\varnothing}x=1.$$

\section{Weighted lattice polynomial functions}

Before introducing the concept of w.l.p.\ function, let us recall the definition of lattice polynomial functions; see for instance Gr\"atzer
\cite[\S I.4]{Grae03}.

\begin{defn}
The class of lattice polynomial functions from $L^n$ to $L$ is defined as follows:
\begin{enumerate}
\item For any $k\in [n]$, the projection $(x_1,\ldots,x_n)\mapsto x_k$ is a lattice polynomial function from $L^n$ to $L$.

\item If $p$ and $q$ are lattice polynomial functions from $L^n$ to $L$, then $p\wedge q$ and $p\vee q$ are lattice polynomial functions from
$L^n$ to $L$.

\item Every lattice polynomial function from $L^n$ to $L$ is constructed by finitely many applications of the rules (1) and (2).
\end{enumerate}
\end{defn}


We now recall that, in a distributive lattice, any lattice polynomial function can be written in disjunctive and conjunctive normal forms, that
is, as a join of meets and dually; see for instance Birkhoff \cite[\S II.5]{Bir67}.

\begin{prop}\label{prop:lp dnf}
Let $p:L^n\to L$ be any lattice polynomial function. Then there are integers $k,l\geqslant 1$ and families $\{A_j\}_{j=1}^k$ and
$\{B_j\}_{j=1}^l$ of nonempty subsets of $[n]$ such that
$$
p(\mathbf{x})=\bigvee_{j=1}^k\; \bigwedge_{i\in A_j}x_i = \bigwedge_{j=1}^l\; \bigvee_{i\in B_j}x_i.
$$
Equivalently, there are nonconstant set functions $\alpha:2^{[n]}\to\{0,1\}$ and $\beta:2^{[n]}\to\{0,1\}$, with $\alpha(\varnothing)=0$ and
$\beta(\varnothing)=1$, such that
$$
p(\mathbf{x})=\bigvee_{\textstyle{S\subseteq [n]\atop \alpha(S)=1}}\bigwedge_{i\in S}x_i=\bigwedge_{\textstyle{S\subseteq
[n]\atop \beta(S)=0}}\bigvee_{i\in S}x_i.
$$
\end{prop}

As mentioned in the introduction, the concept of lattice polynomial function can be generalized by fixing some variables as parameters. Based on
this observation, we naturally introduce the class of w.l.p.\ functions as follows.

\begin{defn}\label{de:wlp}
The class of w.l.p.\ functions from $L^n$ to $L$ is defined as follows:
\begin{enumerate}
\item For any $k\in [n]$ and any $c\in L$, the projection $(x_1,\ldots,x_n)\mapsto x_k$ and the constant function $(x_1,\ldots,x_n)\mapsto c$
are w.l.p.\ functions from $L^n$ to $L$.

\item If $p$ and $q$ are w.l.p.\ functions from $L^n$ to $L$, then $p\wedge q$ and $p\vee q$ are w.l.p.\ functions from $L^n$ to $L$.

\item Every w.l.p.\ function from $L^n$ to $L$ is constructed by finitely many applications of the rules (1) and (2).
\end{enumerate}
\end{defn}

\begin{rem}
Thus defined, w.l.p.\ functions are simply, in the universal algebra terminology, those functions which are definable by polynomial expressions;
see for instance Kaarli and Pixley \cite{KaaPix01} and Lausch and N\"obauer \cite{LauNob73}. Furthermore, these functions are clearly
nondecreasing in each variable.
\end{rem}

Using Proposition~\ref{prop:lp dnf}, we can easily see that any w.l.p.\ function can be written in disjunctive and conjunctive normal forms (see
also Lausch and N\"obauer \cite{LauNob73} and Ovchinnikov \cite{Ovc99}).

\begin{prop}\label{prop:wlp dnf}
Let $p:L^n\to L$ be any w.l.p.\ function. Then there are integers $k,l\geqslant 1$, parameters $a_1,\ldots,a_k,b_1,\ldots,b_l\in L$, and
families $\{A_j\}_{j=1}^k$ and $\{B_j\}_{j=1}^l$ of subsets of $[n]$ such that
$$
p(\mathbf{x})=\bigvee_{j=1}^k\Big(a_j\wedge\bigwedge_{i\in A_j}x_i\Big)=\bigwedge_{j=1}^l\Big(b_j\vee\bigvee_{i\in
B_j}x_i\Big).
$$
Equivalently, there exist set functions $\alpha:2^{[n]}\to L$ and $\beta:2^{[n]}\to L$ such that
$$
p(\mathbf{x})=\bigvee_{S\subseteq [n]}\Big(\alpha(S)\wedge\bigwedge_{i\in S}x_i\Big)=\bigwedge_{S\subseteq
[n]}\Big(\beta(S)\vee\bigvee_{i\in S}x_i\Big).
$$
\end{prop}

It follows from Proposition~\ref{prop:wlp dnf} that any $n$-ary w.l.p.\ function is entirely determined by $2^n$ parameters.

\begin{rem}
Proposition~\ref{prop:wlp dnf} naturally includes the lattice polynomial functions. To see this, it suffices to consider nonconstant set
functions $\alpha:2^{[n]}\to\{0,1\}$ and $\beta:2^{[n]}\to\{0,1\}$, with $\alpha(\varnothing)=0$ and $\beta(\varnothing)=1$.
\end{rem}

\section{Disjunctive and conjunctive normal forms}

We now investigate the link between a given w.l.p.\ function and the parameters that define it.

Let us denote by $p_{\alpha}^{\vee}$ (resp.\ $p_{\beta}^{\wedge}$) the w.l.p.\ function disjunctively (resp.\ conjunctively) defined by the set
function $\alpha:2^{[n]}\to L$ (resp.\ $\beta:2^{[n]}\to L$), that is,
\begin{eqnarray*}
p_{\alpha}^{\vee}(\mathbf{x})&:=&\bigvee_{S\subseteq [n]}\Big(\alpha(S)\wedge\bigwedge_{i\in
S}x_i\Big),\\
p_{\beta}^{\wedge}(\mathbf{x})&:=&\bigwedge_{S\subseteq [n]}\Big(\beta(S)\vee\bigvee_{i\in S}x_i\Big).
\end{eqnarray*}

Of course, the set functions $\alpha$ and $\beta$ are not uniquely determined. For instance, both expressions $x_1\vee(x_1\wedge x_2)$ and $x_1$
represent the same lattice polynomial function.

For any w.l.p.\ function $p:L^n\to L$, define the set functions $\alpha_p:2^{[n]}\to L$ and $\beta_p:2^{[n]}\to L$ as
$\alpha_p(S):=p(\mathbf{e}_S)$ and $\beta_p(S):=p(\mathbf{e}_{[n]\setminus S})$ for all $S\in [n]$. Since $p$ is nondecreasing, $\alpha_p$ is
isotone and $\beta_p$ is antitone.

\begin{lem}\label{lemma:pp}
For any w.l.p.\ function $p:L^n\to L$ we have $p=p^{\vee}_{\alpha_p}=p^{\wedge}_{\beta_p}$.
\end{lem}

\begin{pf*}{Proof.}
Let us establish the first equality. The other one can be proved similarly.

By Proposition~\ref{prop:wlp dnf}, there exists a set function $\alpha:2^{[n]}\to L$ such that $p=p^{\vee}_{\alpha}$. It follows that
$$\alpha_p(T)=\bigvee_{S\subseteq T}\alpha(S)\qquad (T\subseteq [n]).$$ Therefore, we have
\begin{eqnarray*}
p^{\vee}_{\alpha_p}(\mathbf{x}) &=& \bigvee_{T\subseteq [n]}\Big(\alpha_p(T)\wedge\bigwedge_{i\in T}x_i\Big)%
~=~ \bigvee_{T\subseteq [n]}\Big(\bigvee_{S\subseteq T}\alpha(S)\wedge\bigwedge_{i\in T}x_i\Big)\\
&=& \bigvee_{T\subseteq [n]}\bigvee_{S\subseteq T}\Big(\alpha(S)\wedge\bigwedge_{i\in T}x_i\Big)%
~=~ \bigvee_{S\subseteq [n]}\bigvee_{T\supseteq S}\Big(\alpha(S)\wedge\bigwedge_{i\in T}x_i\Big)\\
&=& \bigvee_{S\subseteq [n]}\Big(\alpha(S)\wedge\bigvee_{T\supseteq S}\bigwedge_{i\in T}x_i\Big)%
~=~ \bigvee_{S\subseteq [n]}\Big(\alpha(S)\wedge\bigwedge_{i\in S}x_i\Big)\\
&=& p(\mathbf{x}).\qed
\end{eqnarray*}
\end{pf*}

It follows from Lemma~\ref{lemma:pp} that any $n$-ary w.l.p.\ function is entirely determined by its restriction to $\{0,1\}^n$.

Assuming that $L$ is a chain (that is, $L$ is totally ordered), we now describe the class of all set functions that disjunctively (or
conjunctively) define a given w.l.p.\ function.

\begin{prop}\label{prop:mainalpha}
Assume that $L$ is a chain. Let $p:L^n\to L$ be any w.l.p.\ function and consider two set functions $\alpha:2^{[n]}\to L$ and $\beta:2^{[n]}\to
L$.
\begin{enumerate}
\item We have $p_{\alpha}^{\vee}=p$ if and only if $\alpha^*_p\leqslant\alpha\leqslant \alpha_p$, where the set function $\alpha^*_p:2^{[n]}\to
L$ is defined as
$$
\alpha^*_p(S)=
\begin{cases}
\alpha_p(S), & \mbox{if $\alpha_p(S)>\alpha_p(S\setminus\{i\})$ for all $i\in S$,}\\
0, & \mbox{otherwise.}
\end{cases}
$$
\item We have $p_{\beta}^{\wedge}=p$ if and only if $\beta_p\leqslant\beta\leqslant \beta^*_p$, where the set function $\beta^*_p:2^{[n]}\to L$
is defined as
$$
\beta^*_p(S)=
\begin{cases}
\beta_p(S), & \mbox{if $\beta_p(S)<\beta_p(S\setminus\{i\})$ for all $i\in S$,}\\
1, & \mbox{otherwise.}
\end{cases}
$$
\end{enumerate}
\end{prop}

\begin{pf*}{Proof.}
Let us prove the first assertion. The other one can be proved similarly.

$(\Rightarrow)$ Assume $p_{\alpha}^{\vee}=p$ and fix $S\subseteq [n]$. On the one hand, we have
$$
0\leqslant\alpha(S)\leqslant\bigvee_{K\subseteq S}\alpha(K)=\alpha_p(S).
$$
On the other hand, if $\alpha_p(S)>\alpha_p(S\setminus\{i\})$ for all $i\in S$, then $\alpha(S)=\alpha_p(S)$. Indeed, otherwise, since $L$ is a
chain, there would exist $K^*\varsubsetneq S$ such that
$$
\alpha_p(S)= \bigvee_{K\subseteq S}\alpha(K) = \alpha(K^*) \leqslant \alpha_p(K^*) < \alpha_p(S),
$$
which is a contradiction.

$(\Leftarrow)$ By Lemma~\ref{lemma:pp}, we have $p=p^{\vee}_{\alpha_p}$. Fix $S\subseteq [n]$ and assume there is $i\in S$ such that
$\alpha_p(S)=\alpha_p(S\setminus\{i\})$. Then
$$
\Big(\alpha_p(S\setminus\{i\})\wedge\bigwedge_{j\in S\setminus\{i\}}x_j\Big)\vee \Big(\alpha_p(S)\wedge\bigwedge_{j\in
S}x_j\Big)=\Big(\alpha_p(S\setminus\{i\})\wedge\bigwedge_{j\in S\setminus\{i\}}x_j\Big)
$$
and hence $\alpha_p(S)$ can be replaced with any lower value without altering $p^{\vee}_{\alpha_p}$. Hence
$p^{\vee}_{\alpha_p}=p^{\vee}_{\alpha}$.\qed
\end{pf*}

\begin{exmp}
Assuming that $L$ is a chain, the possible disjunctive expressions of $x_1\vee(x_1\wedge x_2)$ as a 2-ary w.l.p.\ function are given by
$$
x_1\vee(c\wedge x_1\wedge x_2)\qquad (c\in L).
$$
For $c=0$, we retrieve $x_1$ and, for $c=1$, we retrieve $x_1\vee(x_1\wedge x_2)$.
\end{exmp}

We note that, from among all the set functions that disjunctively (or conjunctively) define a given w.l.p.\ function $p$, only $\alpha_p$
(resp.\ $\beta_p$) is isotone (resp.\ antitone). Indeed, suppose for instance that $\alpha$ is isotone. Then, for any $S\subseteq [n]$, we have
$$
\alpha(S)=\bigvee_{K\subseteq S}\alpha(K)=\alpha_p(S),
$$
that is, $\alpha=\alpha_p$.

\section{The discrete Sugeno integral}

Certain w.l.p.\ functions have been considered in the area of nonlinear aggregation and integration. The best known instances are given by the
discrete \emph{Sugeno integral}, which is a particular discrete integration with respect to a \emph{fuzzy measure} (see Sugeno
\cite{Sug74,Sug77}). For a recent survey on the discrete Sugeno integral, see Dubois et al.\ \cite{DubMarPraRouSab01}.

In this section we show the relationship between the discrete Sugeno integral and the w.l.p.\ functions. To this end, we introduce the Sugeno
integral as a function from $L^n$ to $L$. Originally defined when $L$ is the real interval $[0,1]$, the Sugeno integral has different equivalent
representations (see Section~\ref{sec:RepThm}). Here we consider its disjunctive normal representation \cite{Sug74}, which enables us to extend
the original definition of the Sugeno integral to the more general case where $L$ is any bounded distributive lattice.

\begin{defn}
An $L$-valued \emph{fuzzy measure}\/ on $[n]$ is an isotone set function $\mu:2^{[n]}\to L$ such that $\mu(\varnothing)=0$ and $\mu([n])=1$.
\end{defn}

\begin{defn}
Let $\mu$ be an $L$-valued fuzzy measure on $[n]$. The \emph{Sugeno integral}\/ of a function $\mathbf{x}:[n]\to L$ with respect to $\mu$ is
defined by
$$
\sug(\mathbf{x}):=\bigvee_{S\subseteq [n]}\Big(\mu(S)\wedge\bigwedge_{i\in S}x_i\Big).
$$
\end{defn}

Surprisingly, it appears immediately that any function $f:L^n\to L$ is an $n$-ary Sugeno integral if and only if it is a w.l.p.\ function
fulfilling $f(\mathbf{e}_{\varnothing})=0$ and $f(\mathbf{e}_{[n]})=1$. Moreover, as the following proposition shows, any w.l.p.\ function can
be easily expressed in terms of a Sugeno integral.

Recall that, when $n$ is odd, $n=2k-1$, the $n$-ary \emph{median}\/ function is defined in any distributive lattice as the following lattice
polynomial function (see for instance Barbut and Monjardet \cite[Chap.\ IV]{BarMon70})
$$
\med(\mathbf{x}) = \bigvee_{\textstyle{S\subseteq [2k-1] \atop |S|=k}} \bigwedge_{i\in S} x_i = \bigwedge_{\textstyle{S\subseteq [2k-1] \atop
|S|=k}} \bigvee_{i\in S} x_i.
$$

\begin{prop}\label{prop:psug}
For any w.l.p.\ function $p:L^n\to L$, there exists a fuzzy measure $\mu:2^{[n]}\to L$ such that
$$
p(\mathbf{x})=\med\big(p(\mathbf{e}_{\varnothing}),\sug(\mathbf{x}),p(\mathbf{e}_{[n]})\big).
$$
\end{prop}

\begin{pf*}{Proof.}
Let $\mu:2^{[n]}\to L$ be the fuzzy measure which coincides with $\alpha_p$ on $2^{[n]}$ except at $\varnothing$ and $[n]$. Then, we have
\begin{eqnarray*}
\lefteqn{\med\big(p(\mathbf{e}_{\varnothing}),\sug(\mathbf{x}),p(\mathbf{e}_{[n]})\big)}\\
&=& \bigg(\alpha_p(\varnothing)\vee \bigvee_{\textstyle{S\subseteq [n] \atop S\neq\varnothing,
S\neq[n]}}\Big(\mu(S)\wedge\bigwedge_{i\in S}x_i\Big)\vee\Big(\bigwedge_{i\in [n]}x_i\Big)\bigg)\wedge \alpha_p([n])\\
&=& \bigvee_{S\subseteq [n]}\Big(\alpha_p(S)\wedge\bigwedge_{i\in S}x_i\Big)\\
&=& p(\mathbf{x}).\qed
\end{eqnarray*}
\end{pf*}

\begin{cor}
Consider a function $f:L^n\to L$. The following assertions are equivalent:
\begin{enumerate}
\item $f$ is a Sugeno integral.%
\item $f$ is an idempotent w.l.p.\ function, i.e., such that $f(x,\ldots,x)=x$ for all $x\in L$.%
\item $f$ is a w.l.p.\ function fulfilling $f(\mathbf{e}_{\varnothing})=0$ and $f(\mathbf{e}_{[n]})=1$.
\end{enumerate}
\end{cor}

\begin{pf*}{Proof.}
$(1)\Rightarrow (2)\Rightarrow (3)$ Trivial.

$(3)\Rightarrow (1)$ Immediate consequence of Proposition~\ref{prop:psug}.\qed
\end{pf*}

\begin{rem}
As the definition of the w.l.p.\ functions almost coincide with that of the Sugeno integral, certain properties of the Sugeno integral can be
applied as-is or in a slightly extended form to the w.l.p.\ functions. For instance, Proposition~\ref{prop:mainalpha} was already known for the
Sugeno integral (see Marichal \cite{Mar00c}).
\end{rem}

\section{A representation theorem}
\label{sec:RepThm}

Combining Proposition~\ref{prop:psug} with the well-known representations of the Sugeno integral, we easily deduce equivalent representations
for the w.l.p.\ functions.

When $L$ is a chain, for any permutation $\sigma$ on $[n]$, we define the subset
$$
\mathcal{O}_\sigma:=\{\mathbf{x}\in L^n\mid x_{\sigma(1)}\leqslant\cdots\leqslant x_{\sigma(n)}\}.
$$

\begin{thm}\label{thm:represwlp}
Let $p:L^n\to L$ be any w.l.p.\ function. For any $\mathbf{x}\in L^n$, we have
$$
p(\mathbf{x}) = \bigvee_{S\subseteq [n]}\Big(\alpha_p(S)\wedge\bigwedge_{i\in S}x_i\Big) = \bigwedge_{S\subseteq [n]}\Big(\alpha_p(N\setminus
S)\vee\bigvee_{i\in S}x_i\Big).
$$
Moreover, assuming that $L$ is a chain, for any permutation $\sigma$ on $[n]$ and any $\mathbf{x}\in\mathcal{O}_{\sigma}$, setting
$S_{\sigma}(i):=\{\sigma(i),\ldots,\sigma(n)\}$ for all $i\in [n]$, we have
\begin{eqnarray*}
p(\mathbf{x})&=&\bigvee_{i=1}^{n+1} \big(\alpha_p(S_{\sigma}(i))\wedge x_{\sigma(i)}\big)%
~=~\bigwedge_{i=0}^n \big(\alpha_p(S_{\sigma}(i+1))\vee x_{\sigma(i)}\big)\\
&=& \med\big(x_1,\ldots,x_n,\alpha_p(S_{\sigma}(1)),\alpha_p(S_{\sigma}(2)),\ldots,\alpha_p(S_{\sigma}(n+1))\big),
\end{eqnarray*}
with the convention that $x_{\sigma(0)}=0$, $x_{\sigma(n+1)}=1$, and $S_{\sigma}(n+1)=\varnothing$.
\end{thm}

\begin{pf*}{Proof.}
The first part has been established in Lemma~\ref{lemma:pp}. The second part follows from Proposition~\ref{prop:psug} and the following
representations of the Sugeno integral. For any $L$-valued fuzzy measure $\mu$ on $[n]$, we have (see for instance \cite{Mar00c})
\begin{eqnarray*}
\sug(\mathbf{x})&=&\bigvee_{i=1}^n \big(\mu(S_{\sigma}(i))\wedge x_{\sigma(i)}\big)%
~=~ \bigwedge_{i=1}^n \big(\mu(S_{\sigma}(i+1))\vee x_{\sigma(i)}\big)\\
&=& \med\big(x_1,\ldots,x_n,\mu(S_{\sigma}(2)),\mu(S_{\sigma}(3)),\ldots,\mu(S_{\sigma}(n))\big).\qed
\end{eqnarray*}
\end{pf*}

\begin{rem}
It follows from Theorem~\ref{thm:represwlp} that, when the order of the coordinates of $\mathbf{x}$ is known, then $p(\mathbf{x})$ is entirely
determined by $(n+1)$ parameters (instead of $2^n$).
\end{rem}

\section{The median based decomposition formula}

Given a function $f:L^n\to L$ and an index $k\in [n]$, we define the functions $f_k^{0}:L^n\to L$ and $f_k^{1}:L^n\to L$ as
\begin{eqnarray*}
f_k^{0}(\mathbf{x}) &=& f(x_1,\ldots,x_{k-1},0,x_{k+1},\ldots,x_n),\\
f_k^{1}(\mathbf{x}) &=& f(x_1,\ldots,x_{k-1},1,x_{k+1},\ldots,x_n).
\end{eqnarray*}

Clearly, if $f$ is a w.l.p.\ function, so are $f_k^{0}$ and $f_k^{1}$.

Now consider the following system of $n$ functional equations, which we will refer to as the \emph{median based decomposition formula}:
\begin{equation}\label{eq:med k}
f(\mathbf{x})=\med\big(f_k^{0}(\mathbf{x}),x_k,f_k^{1}(\mathbf{x})\big)\qquad (k\in [n])
\end{equation}
This functional system expresses that, for any index $k$, the variable $x_k$ can be totally isolated in $f(\mathbf{x})$ by means of a median
calculated over the variable $x_k$ and the two functions $f_k^{0}$ and $f_k^{1}$, which are independent of $x_k$.

In this final section we establish that this system characterizes the $n$-ary w.l.p.\ functions.

\begin{thm}\label{thm:wlp s1}
The solutions of the median based decomposition formula (\ref{eq:med k}) are exactly the $n$-ary w.l.p.\ functions.
\end{thm}

\begin{pf*}{Proof.}
Recall that the $i$th variable $(i\in [n])$ of a function $f:L^n\to L$ is said to be \emph{effective}\/ if there are two $n$-vectors in $L^n$,
differing only in the $i$th component, on which $f$ takes on different values.

The proof that every function $f:L^n\to L$ satisfying system (\ref{eq:med k}) is a w.l.p.\ function is done by induction on the number of
effective variables of $f$. If $f$ has a single effective variable $x_k$ then, using the $k$th equation of (\ref{eq:med k}), we immediately see
that $f$ is a w.l.p.\ function. The inductive step in then based on the straightforward fact that if $f$ satisfies (\ref{eq:med k}) then, for
any $i\in [n]$, the functions $f_i^0$ and $f_i^1$ also satisfy (\ref{eq:med k}).

Let us now show that any w.l.p.\ function $p:L^n\to L$ fulfills system (\ref{eq:med k}). Let $\mathcal{P}_n$ be the set of nondecreasing
functions $f:L^n\to L$ fulfilling (\ref{eq:med k}). Clearly, $\mathcal{P}_n$ contains all the projection and constant functions from $L^n$ to
$L$. Moreover, we can readily see that if $f,g\in\mathcal{P}_n$ then $f\wedge g\in\mathcal{P}_n$ and $f\vee g\in\mathcal{P}_n$. It follows that
$\mathcal{P}_n$ contains all the w.l.p.\ functions from $L^n$ to $L$.\qed
\end{pf*}

\begin{cor}
For any w.l.p.\ function $p:L^n\to L$ and any $k\in [n]$, we have
$$
p(x_1,\ldots,x_{k-1},p(\mathbf{x}),x_{k+1},\ldots,x_n)=p(\mathbf{x}).
$$
\end{cor}

\begin{pf*}{Proof.}
Using Theorem~\ref{thm:wlp s1} and the fact that $p$ is nondecreasing, we immediately obtain
$$
p(x_1,\ldots,x_{k-1},p(\mathbf{x}),x_{k+1},\ldots,x_n)=\med\big(p_k^{0}(\mathbf{x}),p(\mathbf{x}),p_k^{1}(\mathbf{x})\big)=p(\mathbf{x}).\qed
$$
\end{pf*}

\section{Conclusion}

We have introduced the concept of weighted lattice polynomial functions, which generalize the lattice polynomial functions by allowing some
variables to be fixed as parameters. We have observed that these functions include the class of discrete Sugeno integrals, which have been
extensively used not only in aggregation function theory but also in fuzzy set theory. Finally, we have provided a median based system of
functional equations that completely characterizes the weighted lattice polynomial functions.

Just as special Sugeno integrals (such as the weighted minima, the weighted maxima, and their ordered versions) have already been investigated
and axiomatized (see Dubois et al.\ \cite{DubMarPraRouSab01}), certain subclasses of weighted lattice polynomial functions deserve to be
identified and investigated in detail. This is a topic for future research.

\section*{Acknowledgments}

The author is indebted to Jean-Pierre Barth\'elemy and Stephan Foldes for their comments during the preparation of this paper.

\end{document}